\documentclass[a4paper,11pt]{amsart}
\usepackage{amssymb}
\usepackage[arrow,matrix]{xy}

\title[Bounds on global dimension of piecewise hereditary categories]
{Bounds on the global dimension of certain piecewise hereditary
categories}
\author{Sefi Ladkani}

\address{Einstein Institute of Mathematics, The Hebrew University of Jerusalem, Jerusalem 91904, Israel}
\email{sefil@math.huji.ac.il}

\DeclareMathOperator{\Ext}{Ext}
\DeclareMathOperator{\Hom}{Hom}
\DeclareMathOperator{\ind}{ind}
\DeclareMathOperator{\gldim}{gl.\!dim}
\DeclareMathOperator{\module}{mod}
\DeclareMathOperator{\projd}{pd}
\DeclareMathOperator{\injd}{id}

\newcommand{\cA}{\mathcal{A}}
\newcommand{\cD}{\mathcal{D}}
\newcommand{\cH}{\mathcal{H}}
\newcommand{\cF}{\mathcal{F}}

\newcommand{\bZ}{\mathbb{Z}}

\newcommand{\eps}{\varepsilon}

\newcommand{\dA}{\cD^b(\cA)}
\newcommand{\dH}{\cD^b(\cH)}
\newcommand{\dX}{\cD^b(kX)}
\newcommand{\dY}{\cD^b(kY)}

\theoremstyle{plain}
\newtheorem{theorem}{Theorem}[section]
\newtheorem*{theorem*}{Theorem}
\newtheorem{lemma}[theorem]{Lemma}

\newtheorem{cor}[theorem]{Corollary}

\theoremstyle{definition}

\newtheorem{example}[theorem]{Example}
\newtheorem*{example*}{Example}
\newtheorem{rem}[theorem]{Remark}

\numberwithin{equation}{section}

\begin{document}

\begin{abstract} 
We give bounds on the global dimension of a finite length, piecewise
hereditary category in terms of quantitative connectivity properties of
its graph of indecomposables.

We use this to show that the global dimension of a finite dimensional,
piecewise hereditary algebra $A$ cannot exceed $3$ if $A$ is an
incidence algebra of a finite poset or more generally, a sincere
algebra. This bound is tight.
\end{abstract}

\maketitle

\section{Introduction}

Let $\cA$ be an abelian category and denote by $\dA$ its bounded
derived category. $\cA$ is called \emph{piecewise hereditary} if there
exist an abelian hereditary category $\cH$ and a triangulated
equivalence $\dA \simeq \dH$. Piecewise hereditary categories of
modules over finite dimensional algebras have been studied in the past,
especially in the context of tilting theory, see \cite{Happel88, HRS96,
HRS88}.

It is known \cite[(1.2)]{HRS96} that if $\cA$ is a finite length,
piecewise hereditary category with $n$ non-isomorphic simple objects,
then its global dimension satisfies $\gldim \cA \leq n$. Moreover, this
bound is almost sharp, as there are examples \cite{KSYZ04} where $\cA$
has $n$ simples and $\gldim \cA = n-1$.

In this note we show how rather simple arguments can yield effective
bounds on the global dimension of such a category $\cA$, in terms of
quantitative connectivity conditions on the graph of its
indecomposables, regardless of the number of simple objects.

Let $G(\cA)$ be the directed graph whose vertices are the isomorphism
classes of indecomposables of $\cA$, where two vertices $Q, Q'$ are
joined by an edge $Q \to Q'$ if $\Hom_{\cA}(Q,Q') \neq 0$.

Let $r \geq 1$ and let $\eps = (\eps_0, \dots, \eps_{r-1})$ be a
sequence in $\{+1,-1\}^r$. An \emph{$\eps$-path} from $Q$ to $Q'$ is a
sequence of vertices $Q_0=Q, Q_1, \dots, Q_r = Q'$ such that $Q_i \to
Q_{i+1}$ in $G(\cA)$ if $\eps_i = +1$ and $Q_{i+1} \to Q_i$ if $\eps_i
= -1$.

For an object $Q$ of $\cA$, let $\projd_{\cA} Q = \sup \{ d \,:
\text{$\Ext^d_{\cA}(Q,Q') \neq 0$ for some $Q'$} \}$ and $\injd_{\cA} Q
= \sup \{ d \,: \text{$\Ext^d_{\cA}(Q',Q) \neq 0$ for some $Q'$} \}$ be
the projective and injective dimensions of $Q$, so that $\gldim \cA =
\sup_Q \projd_{\cA} Q$.

\begin{theorem} \label{t:gldim}
Let $\cA$ be a finite length, piecewise hereditary category. Assume
that there exist $r \geq 1$, $\eps \in \{1,-1\}^r$ and an
indecomposable $Q_0$ such that for any indecomposable $Q$ there exists
an $\eps$-path from $Q_0$ to $Q$.

Then $\gldim \cA \leq r+1$ and $\projd_{\cA} Q + \injd_{\cA} Q \leq
r+2$ for any indecomposable $Q$.
\end{theorem}

We give two applications of this result for finite dimensional
algebras.

Let $A$ be a finite dimensional algebra over a field $k$, and denote by
$\module A$ the category of finite dimensional right $A$-modules.
Recall that a module $M$ in $\module A$ is \emph{sincere} if all the
simple modules occur as composition factors of $M$. The algebra $A$ is
called sincere if there exists a sincere indecomposable module.

\begin{cor} \label{c:sincere}
Let $A$ be a finite dimensional, piecewise hereditary, sincere algebra.
Then $\gldim A \leq 3$ and $\projd Q + \injd Q \leq 4$ for any
indecomposable module $Q$ in $\module A$.
\end{cor}

Let $X$ be a finite partially ordered set (\emph{poset}) and let $k$ be
a field. The \emph{incidence algebra} $kX$ is the $k$-algebra spanned
by the elements $e_{xy}$ for the pairs $x \leq y$ in $X$, with the
multiplication defined by setting $e_{xy}e_{y'z} = e_{xz}$ when $y=y'$
and zero otherwise.

\begin{cor} \label{c:incidence}
Let $X$ be a finite poset. If the incidence algebra $kX$ is piecewise
hereditary, then $\gldim kX \leq 3$ and $\projd Q + \injd Q \leq 4$ for
any indecomposable $kX$-module $Q$.
\end{cor}

The bounds in Corollaries~\ref{c:sincere} and~\ref{c:incidence} are
sharp, see Examples~\ref{ex:gldim3} and~\ref{ex:pdid4}.

The paper is organized as follows. In Section~\ref{sec:proofs} we give
the proofs of the above results. Examples demonstrating various aspects
of these results are given in Section~\ref{sec:examples}.

\section{The proofs}
\label{sec:proofs}

\subsection{Preliminaries}

Let $\cA$ be an abelian category. If $X$ is an object of $\cA$, denote
by $X[n]$ the complex in $\dA$ with $X$ at position $-n$ and $0$
elsewhere. Denote by $\ind \cA$, $\ind \dA$ the sets of isomorphism
classes of indecomposable objects of $\cA$ and $\dA$, respectively. The
map $X \mapsto X[0]$ is a fully faithful functor $\cA \to \dA$ which
induces an embedding $\ind \cA \hookrightarrow \ind \dA$.

Assume that there exists a triangulated equivalence $F: \dA \to \dH$
with $\cH$ hereditary. Then $F$ induces a bijection $\ind \dA \simeq
\ind \dH$, and we denote by $\varphi_F: \ind \cA \to \ind \cH \times
\bZ$ the composition
\[
\ind \cA \hookrightarrow \ind \dA \xrightarrow{\sim} \ind \dH = \ind
\cH \times \bZ
\]
where the last equality follows from~\cite[(2.5)]{Keller07}.

If $Q$ is an indecomposable of $\cA$, write $\varphi_F(Q) = (f_F(Q),
n_F(Q))$ where $f_F(Q) \in \ind \cH$ and $n_F(Q) \in \bZ$, so that
$F(Q[0]) \simeq f_F(Q)[n_F(Q)]$ in $\dH$. From now on we fix the
equivalence $F$, and omit the subscript $F$.

\begin{lemma}
The map $f : \ind \cA \to \ind \cH$ is one-to-one.
\end{lemma}
\begin{proof}
If $Q, Q'$ are two indecomposables of $\cA$ such that $f(Q), f(Q')$ are
isomorphic in $\cH$, then $Q[n(Q')-n(Q)] \simeq Q'[0]$ in $\dA$, hence
$n(Q) = n(Q')$, and $Q \simeq Q'$ in $\cA$.
\end{proof}

As a corollary, note that if $A$ and $H$ are two finite dimensional
algebras such that $\cD^b(\module A) \simeq \cD^b(\module H)$ and $H$
is hereditary, then the representation type of $H$ dominates that of
$A$.

We recall the following three results, which were introduced in
\cite[(IV,1)]{Happel88} when $\cH$ is the category of representations
of a quiver.

\begin{lemma}
\label{l:ext} Let $Q, Q'$ be two indecomposables of $\cA$, Then
\[
\Ext^i_{\cA}(Q, Q') \simeq \Ext^{i+n(Q')-n(Q)}_{\cH}(f(Q), f(Q'))
\]
\end{lemma}

\begin{cor}
\label{c:diff} Let $Q, Q'$ be two indecomposables of $\cA$ with
$\Hom_{\cA}(Q,Q') \neq 0$. Then $n(Q') - n(Q) \in \{0,1\}$.
\end{cor}

\begin{lemma}
\label{l:pdid} Assume that $\cA$ is of finite length and there exist
integers $n_0, d$ such that $n_0 \leq n(P) < n_0 + d$ for every
indecomposable $P$ of $\cA$.

If $Q$ is indecomposable, then $\projd_{\cA} Q \leq n(Q)-n_0+1$ and
$\injd_{\cA} Q \leq n_0+d-n(Q)$. In particular, $\gldim \cA \leq d$.
\end{lemma}
\begin{proof}
See~\cite[IV, p.158]{Happel88} or~\cite[(1.2)]{HRS96}.
\end{proof}

\subsection{Proof of Theorem~\ref{t:gldim}}
Let $r \geq 1$, $\eps = (\eps_0,\dots,\eps_{r-1})$ and $Q_0$ be as in
the Theorem. Denote by $r_{+}$ the number of positive $\eps_i$, and by
$r_{-}$ the number of negative ones. Let $F : \dA \to \dH$ be a
triangulated equivalence and write $f = f_F$, $n = n_F$.

Let $Q$ be any indecomposable of $\cA$. By assumption, there exists an
$\eps$-path $Q_0, Q_1, \dots, Q_r = Q$, so by Corollary~\ref{c:diff},
$n(Q_{i+1}) - n(Q_i) \in \{0, \eps_i\}$ for all $0 \leq i < r$. It
follows that $n(Q) - n(Q_0) = \sum_{i=0}^{r-1} \alpha_i \eps_i$ for
some $\alpha_i \in \{0,1\}$, hence
\[
n(Q_0) - r_{-} \leq n(Q) \leq n(Q_0) + r_{+}
\]
and the result follows from Lemma~\ref{l:pdid} with $d=r+1$ and $n_0 =
n(Q_0) - r_{-}$.

\subsection{Variations and comments}

\begin{rem} \label{r:simple}
The assumption in Theorem~\ref{t:gldim} that any indecomposable $Q$ is
the end of an $\eps$-path from $Q_0$ can replaced by the weaker
assumption that any \emph{simple} object is the end of such a path.
\end{rem}
\begin{proof}
Assume that $\eps_{r-1}=1$ and let $Q$ be indecomposable. Since $Q$ has
finite length, we can find a simple object $S$ with $g : S
\hookrightarrow Q$. Let $Q_0, Q_1, \dots, Q_{r-1}, S$ be an $\eps$-path
from $Q_0$ to $S$ with $f_{r-1} : Q_{r-1} \twoheadrightarrow S$.
Replacing $S$ by $Q$ and $f_{r-1}$ by $gf_{r-1} \neq 0$ gives an
$\eps$-path from $Q_0$ to $Q$.

The case $\eps_{r-1} = -1$ is similar.
\end{proof}

\begin{rem} \label{r:diam}
Let $\widetilde{G}(\cA)$ be the undirected graph obtained from $G(\cA)$
by forgetting the directions of the edges. The \emph{distance} between
two indecomposables $Q$ and $Q'$, denoted $d(Q,Q')$, is defined as the
length of the shortest path in $\widetilde{G}(\cA)$ between them (or
$+\infty$ if there is no such path).

The same proof gives that $|n(Q) - n(Q')| \leq d(Q,Q')$ for any two
indecomposables $Q$ and $Q'$. Let $d = \sup_{Q, Q'} d(Q,Q')$ be the
\emph{diameter} of $\widetilde{G}(\cA)$. When $d < \infty$, $\inf_Q
n(Q)$ and $\sup_Q n(Q)$ are finite, and by Lemma~\ref{l:pdid} $\gldim
\cA \leq d+1$ and $\projd_{\cA} Q + \injd_{\cA} Q \leq d+2$ for any
indecomposable $Q$.
\end{rem}

\begin{rem} \label{r:dirsum}
The conclusion of Theorem~\ref{t:gldim} (or Remark~\ref{r:diam}) is
still true under the slightly weaker assumption that $\cA$ is a finite
length, piecewise hereditary category and $\cA = \oplus_{i=1}^{r}
\cA_i$ is a direct sum of abelian full subcategories such that each
graph $G(\cA_i)$ satisfies the corresponding connectivity condition.
\end{rem}

\subsection{Proof of Corollary~\ref{c:sincere}}

Let $A$ be sincere, and let $S_1,\dots,S_n$ be the representatives of
the isomorphism classes of simple modules in $\module A$. Let
$P_1,\dots,P_n$ be the corresponding indecomposable projectives and
finally let $M$ be an indecomposable, sincere module.

Take $r=2$ and $\eps = (-1,+1)$. Now observe that any simple $S_i$ is
the end of an $\eps$-path from $M$, as we have a path of nonzero
morphisms $M \leftarrow P_i \twoheadrightarrow S_i$ since $M$ is
sincere. The result now follows by Theorem~\ref{t:gldim} and
Remark~\ref{r:simple}.

\subsection{Proof of Corollary~\ref{c:incidence}}

Let $X$ be a poset and $k$ a field. A \emph{$k$-diagram} $\cF$ is the
data consisting of finite dimensional $k$-vector spaces $\cF(x)$ for $x
\in X$, together with linear transformations $r_{xx'} : \cF(x) \to
\cF(x')$ for all $x \leq x'$, satisfying the conditions $r_{xx} =
1_{\cF(x)}$ and $r_{xx''} = r_{x'x''}r_{xx'}$ for all $x \leq x' \leq
x''$.

The category of finite dimensional right modules over $kX$ can be
identified with the category of $k$-diagrams over $X$,
see~\cite{Ladkani08}. A complete set of representatives of isomorphism
classes of simple modules over $kX$ is given by the diagrams $S_x$ for
$x \in X$, defined by
\[
S_x(y) = \begin{cases} k & \text{if $y = x$} \\
0 & \text{otherwise}
\end{cases}
\]
with $r_{yy'}=0$ for all $y < y'$. A module $\cF$ is sincere if and
only if $\cF(x) \neq 0$ for all $x \in X$.

The poset $X$ is \emph{connected} if for any $x,y \in X$ there exists a
sequence $x=x_0,x_1,\dots,x_n=y$ such that for all $0 \leq i < n$
either $x_i \leq x_{i+1}$ or $x_i \geq x_{i+1}$.

\begin{lemma} \label{l:kXsincere}
If $X$ is connected then the incidence algebra $kX$ is sincere.
\end{lemma}

\begin{proof}
Let $k_X$ be the diagram defined by $k_X(x)=k$ for all $x \in X$ and
$r_{xx'}=1_k$ for all $x \leq x'$. Obviously $k_X$ is sincere.
Moreover, $k_X$ is indecomposable by a standard connectivity argument;
if $k_X = \cF \oplus \cF'$, write $V = \left\{x \in X \,:\, \cF(x) \neq
0 \right\}$ and assume that $V$ not empty. If $x \in V$ and $x<y$, then
$y \in V$, otherwise we would get a zero map $k \oplus 0 \to 0 \oplus
k$ and not an identity map. Similarly, if $y < x$ then $y \in V$. By
connectivity, $V=X$ and $\cF = k_X$.
\end{proof}

If $X$ is connected, Corollary~\ref{c:incidence} now follows from
Corollary~\ref{c:sincere} and Lemma~\ref{l:kXsincere}. For general $X$,
observe that if $\{X_i\}_{i=1}^{r}$ are the connected components of
$X$, then the category $\module kX$ decomposes as the direct sum of the
categories $\module kX_i$, and the result follows from
Remark~\ref{r:dirsum}.

\begin{cor}
\label{cor:gldim3} Let $X$ and $Y$ be posets such that $\dX \simeq \dY$
and $\gldim kY > 3$. Then $kX$ is not piecewise hereditary.
\end{cor}

\section{Examples}
\label{sec:examples}

We give a few examples that demonstrate various aspects of global
dimensions of piecewise hereditary algebras. In these examples, $k$
denotes a field and all posets are represented by their Hasse diagrams.

\begin{example}[\cite{KSYZ04}]
\label{ex:An} Let $n \geq 2$, $Q^{(n)}$ the quiver
\[
0 \xrightarrow{\alpha_1} 1 \xrightarrow{\alpha_2} 2
\xrightarrow{\alpha_3} \dots \xrightarrow{\alpha_n} n
\]
and $I^{(n)}$ be the ideal (in the path algebra $k Q^{(n)}$) generated
by the paths $\alpha_i \alpha_{i+1}$ for $1 \leq i < n$. By \cite[(IV,
6.7)]{Happel88}, the algebra $A^{(n)} = k Q^{(n)} / I^{(n)}$ is
piecewise hereditary of Dynkin type $A_{n+1}$.

For a vertex $0 \leq i \leq n$, let $S_i$, $P_i$, $I_i$ be the simple,
indecomposable projective and indecomposable injective corresponding to
$i$. Then one has $P_n = S_n$, $I_0 = S_0$ and for $0 \leq i < n$, $P_i
= I_{i+1}$ with a short exact sequence $0 \to S_{i+1} \to P_i \to S_i
\to 0$.

The graph $G(\module A^{(n)})$ is shown below (ignoring the self loops
around each vertex).
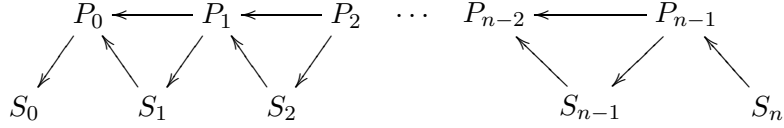
\[
\xymatrix@=0.5pc{ & P_0 \ar[ddl] & & P_1 \ar[ll] \ar[ddl] & & P_2
\ar[ddl] \ar[ll] & {\ldots} & P_{n-2} & & P_{n-1} \ar[ll] \ar[ddl]
\\
\\
S_0 & & S_1 \ar[uul] & & S_2 \ar[uul] & & & & S_{n-1} \ar[uul] & & S_n
\ar[uul] }
\]

Regarding dimensions, we have $\projd S_i = n-i$, $\injd S_i = i$ for
$0 \leq i \leq n$, and $\projd P_i = \injd P_i = 0$ for $0 \leq i < n$,
so that $\gldim A^{(n)} = n$ and $\projd Q + \injd Q \leq n$ for every
indecomposable $Q$. The diameter of $\widetilde{G}(\module A^{(n)})$ is
$n+1$.
\end{example}

The following two examples show that the bounds given in
Corollary~\ref{c:incidence} are sharp.

\begin{example}
\label{ex:gldim3} \emph{A poset $X$ with $kX$ piecewise hereditary and
$\gldim kX = 3$.}

Let $X, Y$ be the two posets:
\[
\begin{array}{ccc}
\xymatrix@=1pc{
& {\bullet} \ar[ddrr] \ar[rr] & & {\bullet} \ar[dr] \\
{\bullet} \ar[dr] \ar[ur] & & & & {\bullet} \\
& {\bullet} \ar[uurr] \ar[rr] & & {\bullet} \ar[ur]
}
& \hspace{40pt} &
\xymatrix@=1pc{
{\bullet} \ar[dr] & & & & {\bullet} \\
& {\bullet} \ar[rr] & & {\bullet} \ar[dr] \ar[ur] \\
{\bullet} \ar[ur] & & & & {\bullet}
}
\\
\\
X & & Y
\end{array}
\]
Then $\dX \simeq \dY$, $\gldim kX = 3$, $\gldim kY = 1$.
\end{example}

\begin{example}
\label{ex:pdid4} \emph{A poset $X$ with $kX$ piecewise hereditary and
an indecomposable $\cF$ such that $\projd_{kX} \cF + \injd_{kX} \cF =
4$.}

Let $X, Y$ be the following two posets:
\[
\begin{array}{ccc}
\xymatrix@=1pc{
& {\bullet} \ar[dr] & & {\bullet} \ar[dr] \\
{\bullet} \ar[dr] \ar[ur] & & {\bullet_x} \ar[dr] \ar[ur] & &
{\bullet} \\
& {\bullet} \ar[ur] & & {\bullet} \ar[ur]
}
& \hspace{40pt} &
\xymatrix@=1pc{
{\bullet} \ar[dr] & & & & {\bullet} \\
& {\bullet} \ar[r] & {\bullet} \ar[r] & {\bullet} \ar[dr] \ar[ur] \\
{\bullet} \ar[ur] & & & & {\bullet}
}
\\
\\
X & & Y
\end{array}
\]
Then $\dX \simeq \dY$, $\gldim kX = 2$, $\gldim kY = 1$ and for the
simple $S_x$ we have $\projd_{kX} S_x = \injd_{kX} S_x = 2$.
\end{example}

We conclude by giving two examples of posets whose incidence algebras
are not piecewise hereditary.

\begin{example}
\emph{A product of two trees whose incidence algebra is not piecewise
hereditary.}

By specifying an orientation $\omega$ on the edges of a (finite) tree
$T$, one gets a finite quiver without oriented cycles whose path
algebra is isomorphic to the incidence algebra of the poset
$X_{T,\omega}$ defined on the set of vertices of $T$ by saying that $x
\leq y$ for two vertices $x$ and $y$ if there is an oriented path from
$x$ to $y$.

A poset of the form $X_{T,\omega}$ is called a \emph{tree}.
Equivalently, a poset is a tree if and only if the underlying graph of
its Hasse diagram is a tree. Obviously, $\gldim k X_{T,\omega}=1$, so
that $k X_{T,\omega}$ is trivially piecewise hereditary. Moreover,
while the poset $X_{T,\omega}$ may depend on the orientation $\omega$
chosen, its derived equivalence class depends only on $T$.

Given two posets $X$ and $Y$, their \emph{product}, denoted $X \times
Y$, is the poset whose underlying set is $X \times Y$ and $(x,y) \leq
(x',y')$ if $x \leq x'$ and $y \leq y'$ where $x,x' \in X$ and $y,y'
\in Y$. It may happen that the incidence algebra of a product of two
trees, although not being hereditary, is piecewise hereditary. Two
notable examples are the product of the Dynkin types $A_2 \times A_2$,
which is piecewise hereditary of type $D_4$, and the product $A_2
\times A_3$ which is piecewise hereditary of type $E_6$.

Consider $X = A_2 \times A_2$ and $Y = D_4$ with the orientations given
below.
\[
\begin{array}{ccc}
\xymatrix@=1pc{
& {\bullet} \ar[dl] \ar[dr] \\
{\bullet} \ar[dr] & & {\bullet} \ar[dl] \\
& {\bullet}
}
& \hspace{40pt} &
\xymatrix@=1pc{
{\bullet} \ar[dr] & & {\bullet} \ar[dl] \\
& {\bullet} \ar[d] \\
& {\bullet}
}
\\
\\
X & & Y
\end{array}
\]
Then $\gldim kX = 2$, $\gldim kY = 1$ and $\dX \simeq \dY$, hence
$\cD^b(k(X \times X)) \simeq \cD^b(k(Y \times Y))$. But $\gldim k(X
\times X) = 4$, so by Corollary~\ref{cor:gldim3}, $Y \times Y$ is a
product of two trees of type $D_4$ whose incidence algebra is not
piecewise hereditary.
\end{example}

\begin{example}
\emph{The converse to Corollary~\ref{c:incidence} is false.}

Let $X$ be the poset
\[
\xymatrix{
{\bullet} \ar[r] \ar[dr] & {\bullet} \ar[r] \ar[dr] & {\bullet} \\
{\bullet} \ar[r] \ar[ur] & {\bullet} \ar[r] \ar[ur] & {\bullet}
}
\]

Then $\gldim kX = 2$, hence $\projd_{kX} \cF \leq 2$, $\injd_{kX} \cF
\leq 2$ for any indecomposable $\cF$, so that $X$ satisfies the
conclusion of Corollary~\ref{c:incidence}. However, $kX$ is not
piecewise hereditary since $\Ext^2_X(k_X, k_X) = k$ does not vanish
(see~\cite[(IV, 1.9)]{Happel88}). Note that $X$ is the smallest poset
whose incidence algebra is not piecewise hereditary.
\end{example}


\end{document}